\documentclass[12pt]{article}

\usepackage{graphicx}
\usepackage{rotating}
\usepackage{amsfonts}
\usepackage{amssymb}
\usepackage{amsmath}
\usepackage{color}
\usepackage{cite}

\def\mt{\mapsto}

\def\hcg{\hat{\cal G}} \def\hch{\hat{\cal H}}
\def\hD{\hat{\Delta}} \def\hpi{\hat{\pi}}
\def\G{\Gamma} 
\def\til{{\tilde I^\L}}
\def\bd{{\bar \D}}
\def\vsh{\vspace{1cm}}
\def\nd{\end{document}}
\def\l{\lambda}
\def\bsq{~~$\blacksquare$}
\def\bbz{Z\!\!\!Z}
\def\({\left(}\def\){\right)}
 
\def\k{\kappa}
\def\rg{\rangle} 
\def\witl{\widetilde{\vert \L\rg}}
\def\bbn{I\!\!N}
\def\r{\rho}
\def\half{{\textstyle{\frac{1}{2}}}}
\def\ha{{\textstyle{\frac{1}{2}}}}
\def\trh{{\textstyle{\frac{3}{2}}}}
\def\frh{{\textstyle{\frac{5}{2}}}}

\def\riga{-\kern-5pt - \kern-5pt -}
\def\white#1{\mathop{\bigcirc}\limits_{#1}}
\def\white#1{\mathop{\circ}\limits_{#1}}
\def\black#1{\mathop{\bullet}\limits_{#1}}

\def\bu{$\bullet~~$}
\def\o{{\bar 0}} \def\I{{\bar 1}}
\def\a{\alpha}
\def\b{\beta}
\def\d{\delta}

\def\D{\Delta} \def\L{\Lambda}
\def\bac{{C\kern-5.5pt I}}
\def\bbc{{C\kern-8pt I}}

\def\eqn{\begin{equation}\label}
\newcommand{\ee}{\end{equation}}

\def\bea{\begin{eqnarray}}

\def\eea{\end{eqnarray}}

\def\bb {\begin {eqnarray}}
\def\eqnn#1{\bb\label{#1}}

\def\nn{\nonumber}

\newcommand{\eqna}[1]{\begin{subequations} \label{#1}
\begin{eqnarray}}
\def\eena{\end{eqnarray}
\end{subequations}}

\input epsf.tex
\newcount\figno
\def\fig#1#2#3{
\par\begingroup\parindent=0pt\leftskip=1cm\rightskip=1cm\parindent=0pt
\baselineskip=11pt \global\advance\figno by 1 
\epsfxsize=#3 \centerline{\epsfbox{#2}} \vskip 12pt
#1\par
\endgroup\par}
\def\figlabel#1{\xdef#1{\the\figno}}
\def\encadremath#1{\vbox{\hrule\hbox{\vrule\kern8pt\vbox{\kern8pt
\hbox{$\displaystyle #1$}\kern8pt} \kern8pt\vrule}\hrule}}

\def\vr{\vert} \def\bbr{I\!\!R}

\def\nt{\noindent}
\def\nl{\hfil\break}

\def\cg{{\cal G}} \def\ch{{\cal H}} 
 \def\ck{{\cal K}} 
  
\def\cp{{\cal P}}

\begin{document}

\begin{center}

{\LARGE {\bf  Positive Energy  Unitary
Irreducible\\[5pt] Representations of the
Superalgebras\\[6pt]  $osp(1|2n,\bbr)$~ and   Character Formulae}}\footnote{To appear in the
Proceedings of the VIII Mathematical Physics Meeting, (Belgrade, 24-31 August 2014).}

\vspace{10mm}

{\bf \large V.K. Dobrev}
\\ \normalsize{Institute of Nuclear Research and Nuclear Energy,\\
 Bulgarian Academy of Sciences,\\
72 Tsarigradsko Chaussee, 1784 Sofia, Bulgaria\\
dobrev@inrne.bas.bg} \vspace{5mm} \\
 {\bf \large I. Salom}
\\ \normalsize{Institute of Physics, University of
Belgrade,\\ Pregrevica 118, 11080 Zemun, Belgrade, Serbia\\
isalom@phy.bg.ac.rs} \vspace{2mm}


\end{center}

\vskip 4mm

\begin{abstract}
We continue the study of positive energy
(lowest weight) unitary irreducible representations of
the superalgebras ~$osp(1|2n,\bbr)$.
We update previous results and present the full list of these UIRs. We
give also some character formulae for these representations.
\end{abstract}


\section{Introduction}

Recently, superconformal field theories in various dimensions are
attracting more interest, in particular, due to their duality to
AdS supergravities.
Until recently only those for ~$D\leq 6$~ were studied since in
these cases the relevant superconformal algebras satisfy \cite{Nahm}
the Haag-Lopuszanski-Sohnius theorem \cite{HLS}.
Thus, such classification was known only for the ~$D=4$~
superconformal algebras ~$su(2,2/N)$ \cite{FF} (for $N=1$),
\cite{DPm,DPu,DPf,DPp} (for arbitrary $N$). More recently, the
classification for ~$D=3$ (for even $N$), $D=5$, and $D=6$ (for
$N=1,2$)~ was given in \cite{Min} (some results are conjectural),
and then the $D=6$ case (for arbitrary $N$) was finalized in \cite{Dosix}.

On the other hand the applications in string theory require the
knowledge of the UIRs of the conformal superalgebras for ~$D>6$.
Most prominent role play the superalgebras $osp(1\vr\,2n)$.
Initially, the superalgebra $osp(1\vr\,32)$
was put forward for $D=10$ \cite{Tow}. Later it was realized that
$osp(1\vr\,2n)$ would fit any dimension, though they are minimal
only for $D=3,9,10,11$ (for $n=2,16,16,32$, resp.) \cite{DFLV}.
In all cases we need to find first the UIRs of
~$osp(1\vr\, 2n,\bbr)$~ which study was started in \cite{DZ} and \cite{DMZZ}.

In the present paper we intend to finalize the unitarity classification of \cite{DZ}
and in addition to provide some character formulae.

Since this paper is a sequel of \cite{DZ}, where there is extensive literature, and for the lack of space we
only update the supersymmetry literature (for $D>2$) after 2004, cf.
\cite{Bei,Sok,BiKo,Var,Raj,Gen,Mil,Met,Ter,Nak,Dol,Do-ch,Ber,RySa,Pol,Bos,HeYa,Sol,Tac,LSV,Die,AnTu,Yon,HoKa,BSZ,GKT,Knu,Ras,Tor,GoSk,Uhl,CGL,HaPe,JuSi,ALS,NeSa,BGH,QuSc,Bui,SpVa,LiSt,BeTs,Cou,MaMo,FSW,DMV,All,GKN}

\section{Representations of the superalgebras $\ osp(1\vr\, 2n)\ $
and $\ osp(1\vr\, 2n,\bbr)$}

\subsection{The setting}

\nt Our basic references for Lie superalgebras are \cite{Kab,Kc},
although in this exposition we follow \cite{DZ}.

The even subalgebra of ~$\cg = osp(1\vr\, 2n,\bbr)$~ is the algebra
~$sp(2n,\bbr)$ with maximal compact subalgebra ~$\ck = u(n) \cong
su(n) \oplus u(1)$.

We label the relevant representations of ~$\cg$~
by the signature:
\eqn{sgn}\chi ~=~ [\, d\,;\,a_1\,,...,a_{n-1}\,] \ee
where ~$d$~ is the conformal weight, and ~$a_1,...,a_{n-1}$~
are non-negative integers which are Dynkin labels of the
finite-dimensional UIRs of the subalgebra $su(n)$ (the simple
part of $\ck$).

In \cite{DZ}  were classified (with some omissions to be spelled out
below) the positive energy (lowest weight) UIRs of ~$\cg$~ following
the methods used for the $D=4,6$ conformal superalgebras, cf.
\cite{DPm,DPu,DPf,DPp,Dosix}, resp. The main tool was an adaptation
of the Shapovalov form \cite{Sha} on the Verma modules ~$V^\chi$~
over the complexification ~$\cg^\bac ~=~ osp(1\vr\, 2n)$~ of ~$\cg$.

\subsection{Root systems}

\nt
We recall some facts about ~$\cg^\bac ~=~ osp(1\vr\, 2n)$
(denoted $B(0,n)$ in \cite{Kab}) as used in \cite{DZ}.
The  root systems are given in terms of
~$\d_1\,\dots,\d_{n}\,$, ~$(\d_i,\d_j) ~=~
\d_{ij}\,$, ~$i,j=1,...,n$.
The even and odd roots systems are \cite{Kab}:
\eqnn{dmn}
\D_\o  &=&  \{ \pm\d_i\pm\d_j\ , ~1\leq i< j\leq n\ ,
~~\pm 2\d_i\ , ~1\leq i \leq n \} ~,
\\ \D_\I  &=&  \{ \pm\d_i\ , ~1\leq i \leq n
\} \nn\eea
(we remind that the signs ~$\pm$~ are not correlated).
We shall use the following distinguished
simple root system \cite{Kab}:
\eqn{ssdm}
 \Pi ~=~ \{\, \d_1-\d_2 \, ,\,
 , \dots, \d_{n-1}- \d_{n} \, ,\, \d_n\, \} \ , \ee
or introducing standard notation for the simple roots:
\eqnn{ssdma}
\Pi ~&=&~ \{\,\a_1\,,...,\,\a_{n} \, \}\ , \\
&&\a_{j} ~=~ \d_{j}-\d_{j+1}\ , \quad
j=1,...,n-1 \ , \quad \a_{n} ~=~ \d_{n}\ . \nn\eea
The root ~$\a_n = \d_n$~ is odd, the other simple roots
are even.
The Dynkin diagram is:
\eqn{dynk}
\vbox{\offinterlineskip\baselineskip=10pt
\halign{\strut#
\hfil
& #\hfil
\cr &\cr
&\ $\white{{1}}
\riga \cdots \riga
\white{{n-1}} =\kern-2pt\Longrightarrow
\black{n}$
\cr }} \ee
The black dot is used to signify that the simple odd root is not
nilpotent, otherwise a gray dot would be used \cite{Kab}. In fact, the
superalgebras ~$B(0,n) = osp(1\vr\, 2n)$~ have no nilpotent generators
unlike all other types of basic classical Lie superalgebras \cite{Kab}.

The corresponding to ~$\Pi$~ positive root system is:
\eqn{psdm}\D_\o^+ ~=~ \{ \d_i\pm\d_j\ , ~1\leq i< j\leq n\ ,
 ~~2\d_i \ , ~1\leq i \leq n \} ~,
\qquad \D_\I^+ ~=~ \{ \d_i \ , ~1\leq i \leq n \} \ee
We record how the elementary functionals are expressed through
the simple roots:
\eqn{back} \d_k ~=~ \a_k + \cdots + \a_n \ .\ee

From the point of view of representation theory more relevant is the restricted root system, such that:
\eqnn{rrr} && \bd^+ ~=~ \bd^+_\o \cup \D_\I^+ \ , \\
&& \bd^+_\o ~\equiv ~ \{ \a\in\D^+_\o ~\vert ~\half\a\notin\D^+_\I
\} ~=~ \{ \d_i\pm\d_j\ , ~1\leq i< j\leq n \} \eea

The superalgebra ~$\cg ~=~ osp(1\vr\, 2n,\bbr)$~ is a split real form
of $osp(1\vr\, 2n)$ and has the same root system.

The above simple root system is also the simple root system of the complex
simple Lie algebra ~$B_n$ (dropping the distinction between even and odd roots) with
Dynkin diagram:
\eqn{dynko}
\vbox{\offinterlineskip\baselineskip=10pt
\halign{\strut#
\hfil
& #\hfil
\cr &\cr
&\ $\white{{1}}
\riga \cdots \riga
\white{{n-1}} =\kern-2pt\Longrightarrow
\white{n}$
\cr }} \ee
Naturally, for the $B_n$ positive root system we drop the roots $2\d_i\,$~
\eqn{psdmb}\D^+_{\rm B_n} ~=~ \{ \d_i\pm\d_j\ , ~1\leq i< j\leq n\ ,
 ~~\d_i \ , ~1\leq i \leq n \} ~\cong ~\bd^+ \ee
This shall be used essentially below.

\subsection{Lowest weight through the signature}

\nt
Besides \eqref{sgn} we shall use the Dynkin-related labelling:
\eqn{rrrr}
(\L , \a_k^\vee ) ~= -\, a_k \ , ~~~ 1\leq k \leq n\ ,   \ee
where ~$\a_k^\vee \equiv 2\a_k/(\a_k,\a_k)$, and
the minus signs
are related to the fact that we work with
lowest weight Verma modules (instead of the highest weight modules
used in \cite{Kc}{}) and to Verma module reducibility
w.r.t. the roots ~$\a_k$ (this is explained in detail in \cite{DPf,DZ}).

Obviously, ~$a_n$~ must be related to  the conformal weight ~$d$~  which is
  a matter of
normalization so as to correspond to some known cases.
Thus, our choice is:
\eqn{cftw} a_n ~=~ - 2 d - a_1 - \cdots - a_{n-1} \ . \ee

 The actual Dynkin labelling is given by:
\eqn{dynks} m_k ~=~ (\r - \L, \a_k^\vee )\ee
where ~$\r\in\ch^*$~ is  given by the
difference of the half-sums $\r_\o\,,\r_\I$
of the even, odd, resp., positive roots (cf. (\ref{psdm}):
\eqnn{hlf} \r ~&\doteq&~ \r_\o - \r_\I
~=~ (n-\half)\d_1 + (n-\trh) \d_2 + \cdots + \trh \d_{n-1} +
\half \d_n \ , \\
&&\r_\o\ =\ n \d_1 + (n-1) \d_2 + \cdots +
2\d_{n-1} + \d_n \ , \nn\\
&&\r_\I \ =\ \half (\d_1 + \cdots + \d_n)\ . \nn\eea

Naturally, the value of $\r$ on the simple roots is 1:
~$(\r,\a^\vee_i)=1$, $i=1,...,n$.

Unlike $a_k\in\bbz_+$ for $k<n$ the value of $a_n$ is arbitrary.
In the cases when $a_n$ is also a non-negative integer, and then ~$m_k\in\bbn$ ($\forall k$) the
corresponding irreps are the finite-dimensional irreps of $\cg$ (and of $B_n$).

Having in hand the values of ~$\L$~ on the basis we can
recover them for any element of ~$\ch^*$.

We shall need only ~$(\L , \b^\vee)$~ for all positive roots $\b$
 as given in \cite{DZ}:
\eqnn{lorr}
(\L , (\d_i-\d_j)^\vee ) ~&=&~ (\L , \d_i-\d_j ) ~=~
- a_i - \cdots - a_{j-1}  \\
(\L , (\d_i+\d_j)^\vee ) ~&=&~ (\L , \d_i+\d_j ) ~=~
2d\ +\ a_1 + \cdots + a_{i-1} - a_j - \cdots - a_{n-1} \qquad
\nn\\
(\L , \d_i^\vee ) ~&=&~ (\L , 2\d_i ) ~=~
2d \ +\ a_1 + \cdots + a_{i-1} - a_i - \cdots - a_{n-1}
\nn\\
(\L , (2\d_i)^\vee ) ~&=&~ (\L , \d_i ) ~=~
d \ +\ \half ( a_1 + \cdots + a_{i-1} - a_i - \cdots - a_{n-1} ) \nn
\eea

\subsection{Verma modules}

\nt
To introduce Verma modules we use the standard
triangular decomposition:
\eqn{trig} \cg^\bac ~=~ \cg^+ \oplus \ch \oplus \cg^- \ee
where $\cg^+$, $\cg^-$, resp., are the subalgebras corresponding
to the positive, negative, roots, resp., and $\ch$ denotes the
Cartan subalgebra.

We consider lowest weight Verma modules,
so that ~$V^\L ~ \cong U(\cg^+) \otimes v_0\,$,
 where ~$U(\cg^+)$~ is the universal enveloping algebra of $\cg^+$,
and ~$v_0$~ is a lowest weight vector $v_0$ such that:
\eqnn{low}
 Z \ v_0\ &=&\ 0 \ , \quad Z\in \cg^- \nn\\
 H \ v_0 \ &=&\ \L(H)\ v_0 \ , \quad H\in \ch \ .\eea
Further, for simplicity we omit the sign
~$\otimes \,$, i.e., we write $p\,v_0\in V^\L$ with $p\in U(\cg^+)$.

Adapting the criterion of \cite{Kc} (which generalizes the
BGG-criterion \cite{BGG} to the super case) to lowest weight
modules,  one finds that a Verma module ~$V^\L$~ is reducible w.r.t.
the positive root ~$\b$~ iff the following holds \cite{DZ}:
\eqn{odr} (\r - \L, \b^\vee) = m_\b\ , \qquad \b \in \Delta^+ \ ,
\quad m_\b\in\bbn \ . \ee

If a condition from (\ref{odr}) is fulfilled then ~$V^\L$~ contains
a submodule which is a Verma module ~$V^{\L'}$~ with shifted weight
given by the pair ~$m,\b$~: ~$\L' ~=~ \L + m\b$. The embedding of
~$V^{\L'}$~ in ~$V^\L$~ is provided by mapping the lowest weight
vector ~$v'_0$~ of ~$V^{\L'}$~ to the ~{\bf singular vector}
~$v_s^{m,\b}$~ in ~$V^\L$~ which is completely determined by the
conditions: \eqnn{lowp}
 X \ v_s^{m,\b}\ &=&\ 0 \ , \quad X\in \cg^- \ , \nn\\
 H \ v_s^{m,\b} \ &=&\ \L'(H)\ v_0 \ , \quad H\in \ch \ ,
~~~\L' ~=~ \L + m\b\ .\eea Explicitly, ~$v_s^{m,\b}$~ is given by a
polynomial in the positive root generators \cite{Dob,DPf}:
\eqn{sing} v_s^{m,\b} ~=~ P^{m,\b} \,v_0 \ , \quad P^{m,\b}\in
U(\cg^+)\ . \ee Thus, the submodule ~$I^\b$~ of ~$V^\L$~ which is
isomorphic to ~$V^{\L'}$~ is given by ~$U(\cg^+)\, P^{m,\b}
\,v_0\,$.

Note that the Casimirs of $\cg^\bac$ take the same values on
~$V^\L$~ and ~$V^{\L'}$.

Certainly, \eqref{odr} may be fulfilled for several positive roots
(even for all of them). Let ~$\D_\L$~ denote the set of all positive
roots for which \eqref{odr}  is fulfilled, and let us denote:
~$\tilde I^\L ~\equiv~ \cup_{\b\in\D_\L} I^{\b}\,$. Clearly,
~$\tilde I^\L$~ is a proper submodule of ~$V^\L$. Let us also denote
~$F^\L ~\equiv V^\L / \tilde I^\L$.

Further we shall use also the following notion. The singular vector
~$v_1$~ is called ~{\bf descendant}~ of the singular vector
~$v_2\notin \bbc v_1$~ if there exists a homogeneous polynomial
~$P_{12}$~ in ~$U(\cg^+)$~ so that ~$v_1 ~=~ P_{12} ~v_2\,$. Clearly,
in this case we have: ~$I^1 ~\subset ~I^2\,$, where ~$I^k$~ is the
submodule generated by $v_k\,$.

The Verma module ~$V^\L$~ contains a unique proper maximal submodule
~$I^\L$ ($\supseteq \til$) \cite{Kc,BGG}. Among the lowest weight
modules with lowest weight ~$\L$~ there is a unique irreducible one,
denoted by ~$L_\L$, i.e., ~$L_\L ~=~ V^\L/I^\L$. (If ~$V^{\L}$~ is
irreducible then ~$L_\L = V^\L$.)

It may happen  that the maximal submodule ~$I^\L$~ coincides with
the submodule ~$\tilde I^\L$~ generated by all singular vectors.
This is, e.g., the case for all Verma modules if ~rank~$\cg ~\leq
2$, or when \eqref{odr} is fulfilled for all simple roots (and, as a
consequence for all positive roots). Here we are interested in the
cases when ~$\tilde I^\L$~ is a proper submodule of ~$I^\L$. We need
the following notion.

{\bf Definition:} \cite{BGG,Docond,Do-KL} ~{\it  Let ~$V^\L$~ be a
reducible Verma module. A vector ~$v_{{\rm ssv}} ~\in ~V^\L$~ is
called a ~{\bf subsingular vector}~ if ~$v_{{\rm su}} ~\notin ~
\tilde I^\L$~ and the following holds:} \eqn{subs} X~v_{{\rm su}}
~~\in ~~ \tilde I^\L \,, \quad \forall X\in\cg^-  \ee

\medskip

Going from the above more general definitions to ~$\cg$~ we recall
that in \cite{DZ}  it was established that from \eqref{odr} follows
that the Verma module ~$V^{\L(\chi)}$~ is reducible if one of the
following relations holds (following the order of (\ref{lorr}):
\eqna{redd}
\bbn \ni m^-_{ij} &=& j-i + a_i + \cdots + a_{j-1} \\
\bbn \ni m^+_{ij} &=& 2n -i-j +1 + a_j + \cdots + a_{n-1}
- a_1 - \cdots - a_{i-1} -2d \qquad \\
\bbn \ni m_i &=& 2n-2i+1 + a_i + \cdots + a_{n-1} - a_1 +
\cdots - a_{i-1} - 2d\  \\
\bbn \ni m_{ii} &=& n-i+\half(1 + a_i + \cdots + a_{n-1} - a_1 +
\cdots - a_{i-1}) - d \ . \eena
Further we shall use  the fact from \cite{DZ}  that we may eliminate the reducibilities and
embeddings related to the roots ~$2\d_i\,$.  Indeed,
since  $m_i ~=~ 2m_{ii}\,$,   whenever (\ref{redd}{d}) is
fulfilled also (\ref{redd}{c}) is fulfilled.

For further use we introduce notation for the root vector ~$X^+_j\,\in\cg^+$, ~$j=1,\ldots,n$,
~corresponding to the simple root ~$\a_j\,$. Naturally, ~$X^-_j\,\in\cg^-$~ corresponds to ~$-\a_j\,$.

Further, we notice that all reducibility conditions in
(\ref{redd}{a}) are fulfilled. In particular, for the simple roots
from those condition (\ref{redd}{a}) is fulfilled with
~$\b\to\a_i=\d_i-\d_{i+1}\,$, $i=1,...,n-1$ and ~$m^-_i ~\equiv~
m^-_{i,i+1} ~=~ 1 + a_i\,$. The corresponding submodules ~$I^{\a_i}
~=~ U(\cg^+)\,v^i_s\,$, where ~$\L_i ~=~ \L + m^-_i \a_i$~ and
~$v^i_s ~=~ (X^+_i)^{1+a_i}\, v_0\,$. These submodules generate an
invariant submodule which we denote by ~$I^\L_c\, \subset \tilde
I^\L$. Since these submodules are nontrivial for all our signatures
in the question of unitarity instead of ~$V^\L$~ we shall consider
also the factor-modules: \eqn{fcc} F_c^\L ~=~ V^\L\, /\, I^\L_c
~\supset ~F^\L\ . \ee We shall denote the lowest weight vector of
~$F_c^\L$~ by ~$\vr \L_c\rg$~ and the singular vectors above become
null conditions in~ $F_c^\L$~: \eqn{nullm} (X^+_i)^{1+a_i}\,
 \vr \L_c\rg ~=~ 0 \ , \quad i=1,...,n-1. \ee

If the Verma module ~$V^\L$~ is not reducible w.r.t. the other
roots, i.e., (\ref{redd}{b,c,d}) are not fulfilled, then
~$F_c^\L=F^\L$~ is irreducible and is isomorphic to the irrep
~$L_\L$~ with this weight.

In fact, for the factor-modules reducibility is controlled by the value of ~$d$,
or in more detail:

The maximal ~$d$~ coming from the different possibilities in
(\ref{redd}{b}) are obtained for $m^+_{ij}=1$ and they are:
\eqn{boun}
d_{ij} ~\equiv~ n + \half(a_j + \cdots + a_{n-1}
- a_1 - \cdots - a_{i-1}-i-j) \ ,
\ee
the corresponding root being ~$\d_i+\d_j\,$.

The maximal ~$d$~ coming from the different possibilities in
(\ref{redd}{c,d}), resp., are obtained for $m_{i}=1$, $m_{ii}=1$,
resp., and they are:
\eqnn{bounz}
&&d_{i} ~\equiv~ n -i + \half(a_i + \cdots + a_{n-1}
- a_1 - \cdots - a_{i-1}) \ ,
\\
&&d_{ii} ~=~ d_i - \half \ ,
\nn\eea
the corresponding roots being ~$\d_i\,,2\d_j\,$, resp.

There are some orderings between these maximal reduction points \cite{DZ}:
\eqnn{compz}
d_1 ~&>&~ d_2 ~>~ \cdots ~>~ d_n \ , \\
d_{i,i+1} ~&>&~ d_{i,i+2} ~>~ \cdots ~>~ d_{in} \ , \nn\\
d_{1,j} ~&>&~ d_{2,j} ~>~ \cdots ~>~ d_{j-1,j} \ , \nn\\
d_i ~&>&~ d_{jk} ~>~ d_\ell
\ , \qquad i\leq j <k \leq \ell \ . \nn\eea

Obviously the first reduction point is:
\eqn{frp} d_{1} ~=~ n -1 + \half(a_1 + \cdots + a_{n-1})  \ .
\ee

\section{Unitarity}

The first results on the unitarity were given in \cite{DZ}. These were not
complete so the statement below should be called~ {\it Dobrev-Zhang-Salom Theorem}.

\nt {\bf Theorem:}~~ All positive energy unitary irreducible
representations of the superalgebras ~$osp(1\vr\, 2n,\bbr)$~ characterized
by the signature ~$\chi$~ in (\ref{sgn}) are obtained for real ~$d$~
and are given as follows:
\eqnn{unto}
&&d ~\geq~ n -1 + \half(a_1 + \cdots + a_{n-1})   ~=~ d_{1}
\ , \quad   a_1\neq 0 \ , \\
&&d ~\geq~  n - \trh + \half(a_2 + \cdots + a_{n-1}  )   ~=~ d_{12}
\ , \quad a_1 = 0,\  a_2\neq 0\ , \nn\\
&&d ~=~  n - 2 + \half(a_2 + \cdots + a_{n-1}  )   ~=~d_{2}  > d_{13}
\ , \quad a_1 = 0,\  a_2\neq 0\ ,  \\
&&d ~\geq~  n - 2 + \half(a_3 + \cdots + a_{n-1}  )   ~=~ d_{2} ~=~ d_{13}
\ , \quad a_1 = a_2=0,\  a_3\neq 0\ , \nn\\
&&d ~=~  n - \frh + \half(a_3 + \cdots + a_{n-1}  )   ~=~ d_{23} ~>~ d_{14}
\ , \quad a_1 = a_2=0,\  a_3\neq 0\ , \nn\\
&&d ~=~  n - 3 + \half(a_3 + \cdots + a_{n-1}  )   ~=~ d_3 ~=~ d_{24} ~>~ d_{15}
\ , \quad a_1 = a_2=0,\  a_3\neq 0\ , \nn\\
&& ... \nn\\
&& ... \nn\\
&&d ~\geq~   n - 1 -\k + \half(a_{2\k+1} + \cdots + a_{n-1})
\ , \quad a_1 = ... = a_{2\k}=0,\  a_{2\k+1}\neq 0\ , \nn\\
&& \qquad\qquad \k = \half, 1, ..., \half(n-1) \ , \nn\\
&&d ~=~   n - \trh - \k  + \half(a_{2\k+1} + \cdots + a_{n-1})
\ , \quad a_1 = ... = a_{2\k}=0,\  a_{2\k+1}\neq 0\ , \nn\\
&& ... \nn\\
&&d ~=~   n - 1-2\k + \half(a_{2\k+1} + \cdots + a_{n-1})
\ , \quad a_1 = ... = a_{2\k}=0,\  a_{2\k+1}\neq 0\ , \nn\\
&& ... \nn\\
&& ... \nn\\
 &&d ~\geq~   \half (n-1)\ , \quad a_1 ~=~ ... ~=~ a_{n-1} ~=~ 0 \nn\\
&&d ~=~   \half (n-2)\ , \quad a_1 ~=~ ... ~=~ a_{n-1} ~=~ 0 \nn\\
&& ... \nn\\
&&d ~=~   \half \ , \quad a_1 ~=~ ... ~=~ a_{n-1} ~=~ 0 \nn\\
&&d ~=~  0 \ , \quad a_1 ~=~ ... ~=~ a_{n-1} ~=~ 0 \nn
 \eea
\nt
{\bf Proof:}~~  The statement of the Theorem for ~$d ~>~ d_1$~ is
clear in \cite{DZ}  from the general considerations since this is the First reduction point.
For ~$d ~=~ d_1$~ (also following \cite{DZ}) we
have the first zero norm state which is naturally given by the
corresponding singular vector ~$v^1_{\d_1} ~=~ \cp^{1,\d_1}\
v_0\,$. In fact, all states of the embedded submodule
~$V^{\L+\d_1}$~ built on ~$v^1_{\d_1}$~ have zero norms.
Due to the above singular vector we have the following
additional null condition in ~$F_c^\L$~:
\eqn{nula} \cp^{1,\d_1}\ \witl ~=~ 0 \ . \ee
The above condition  factorizes the submodule built on
~$v^1_{\d_1}\,$.
There are no other vectors with zero norm at $d=d_1$ since by a
general result \cite{Kc}, the elementary embeddings between Verma
modules are one-dimensional.
Thus, ~$F^\L$~ is the UIR ~$L_\L ~=~ F^\L$.\nl
Below ~$d<d_1$~ there is no unitarity for ~$a_1\neq 0$.
On the other hand (as shown in \cite{DZ})
for $a_1= 0$ the singular vector $v^1_{\d_1}$ is
descendant of the compact root singular vector  ~$X_1^+\,v_0$~
which is already factored out for ~$a_1 = 0$.
Thus, below we set ~$a_1  = 0$.

The next reducibility point is ~$d=d_{12}= n -\trh +\half(a_2+\cdots+a_{n-1})$.
The corresponding root is  ~$\d_1+\d_{2} ~=~
\a_1 + 2\a_{2} +\cdots + 2\a_n\,$.
The corresponding singular vector is ~$v^1_{\d_1+\d_{2}}
~=~ \cp^{1,\d_1+\d_{2}}\ v_0\,$.
All states of the embedded submodule
~$V^{\L+\d_1+\d_{2}}$~ built on ~$v^1_{\d_1+\d_{2}}$~ have
zero norms for ~$d ~=~ d_{12}\,$.
Due to the above singular vector we have the following
additional null condition in ~$F_c^\L$~:
\eqn{nuli} \cp^{1,\d_1+\d_{2}}\ \witl ~=~ 0 \ ,
\qquad d ~=~d_{12}\, . \ee
The above conditions factorizes the submodule built on
~$v^1_{\d_1+\d_{2}}\,$.
Thus, ~$F_c^\L$~ is the UIR ~$L_\L ~=~ F_c^\L$.

Below ~$d<d_{12}$~ there is no unitarity for ~$a_2\neq 0$,
except at the isolated point:
~$d_2= n-2 +\half(a_2+\cdots+a_{n-1})$.
At the latter point there is a singular vector ~$v^1_{\d_{2}}$~
which must be factored for unitarity. In addition, the previous
singular vector is descendant of ~$v^1_{\d_{2}}$~
and the compact root singular vector   ~$X_1^+\,v_0$.

Further, for
for $a_2= 0$ the singular vectors $v^1_{\d_1+\d_{2}}$ and ~$v^1_{\d_{2}}$~ are
descendants of the compact root singular vectors   ~$X_1^+\,v_0$~
and ~$X_2^+\,v_0$~
which are factored out for ~$a_1 = a_2= 0$.
Thus, below we set also ~$a_2  = 0$~ and there would be no obstacles
for unitarity until the next reducibility points (coinciding due $a_2  = 0$):
~$d_2=d_{13}= n-2 +\half(a_3+\cdots+a_{n-1})$.
 The singular vector for ~$d=d_{13}  $~ and ~$m=1$~ has weight
~$\d_1+\d_3=\a_1+\a_2+2\a_3 + \cdots 2\a_n$~ and
 for ~$a_1=0$~ it is a descendant of the compact root singular vector ~$X_1\, v_0\,$\ \cite{Dos}.
 However, at ~$d_2=d_{13}$~ there is a subsingular vector which must   be factored for
unitarity.
 For ~$  d <d_2=d_{13}$~ and ~$a_3\neq 0$~ the norm of that subsingular vector
 is  negative, and there will be no unitarity   except at some  lower reducibility points.

For ~$d_{23} = n - \frac{5}{2} +\half(a_3+\cdots+a_{n-1})$~ there is singular vector
~$v^1_{\d_{2}+\d_3}$~ of weight ~$\d_2+\d_3=\a_2+2\a_3 + \cdots 2\a_n$~ \cite{Dos} which must be factored for unitarity.
The previous subsingular vector is also factored out since it is descendant of ~$v^1_{\d_{2}+\d_3}$~ and
compact root singular vectors.

Further on, the Proof goes on similar lines. We list the points at which there are
subsingular vectors - these happen when reducibility points coincide due the
zero values of some ~$a_i$~:
\eqnn{subsi}
&&d_2=d_{13}= n-2 +\half(a_3+\cdots+a_{n-1}) \ , \quad a_1=a_2 =0\ , \\
&&d_{23}=d_{14}= n-5/2 +\half(a_4+\cdots+a_{n-1}) \ , \quad a_1=a_2 =a_3=0, \ n>3, \nn\\
&&d_3=d_{24}= d_{15}= n-3 +\half(a_5+\cdots+a_{n-1}) \ , \quad a_1=a_2 =a_3 =a_4=0 ,\nn\\
&&\qquad\qquad\qquad\qquad \qquad\qquad\qquad\qquad n>3, \nn\\
&& ... \nn\\
&&d_j=d_{1,2j-1}= d_{2,2j-2}= \cdots = d_{j-1,j+1}=
n-j +\half(a_{2j-1}+\cdots+a_{n-1}) \ , \nn\\ && \qquad\qquad\qquad
a_1=\cdots = a_{2j-2} =0 , \quad j<n , \nn\\
&&d_{j,j+1}=d_{1,2j}= d_{2,2j-1}= \cdots = d_{j-1,j+2}= n-j-\half +
\half(a_{2j} +\cdots+a_{n-1}) \ , \nn\\ && \qquad\qquad\qquad
 a_1=\cdots = a_{2j-1} =0 , \quad j<n-1. \nn
\eea
Above it is understood that ~$a_j\equiv 0$~ for ~$j\geq n$.\nl
At the points of the subsingular vectors the associated singular vectors
are factored out automatically. This happens also when the
subsingular vectors are inside a continuous part of the unitarity spectrum.\bsq

The Proof above is not as explicit as we would like it to be, but due to the lack of space
we postpone it to \cite{DoSa}. Below we give separately and explicitly the case ~$n=3$.

\bigskip

{\bf Example: n=3.}~~ For ~$n=3$~ f-la (\ref{compz}) simplifies to:
\eqn{fig3} \includegraphics[bb=-75 -10 1000 40]{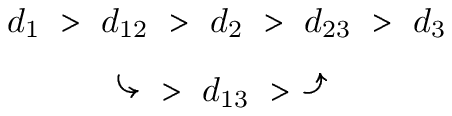} \ee

The Theorem  now reads:\nl
\eqnn{untr}
&&d ~\geq~  2 + \half(a_1 + a_{2})   ~=~d_{1}
\ , \quad   a_1\neq 0 \ , \\
&&d ~\geq~ \trh +  \half a_2   ~=~d_{12}
\ , \quad a_1 = 0,\  a_2\neq 0\ , \nn\\
&&d ~=~  1 + \half a_2   ~=~ d_{2} > d_{13}
\ , \quad a_1 = 0,\  a_2\neq 0\ , \nn\\
&&d ~\geq~ 1 ~=~ d_2 = d_{13}  \ , \quad a_1 = a_2= 0\ , \nn\\
 &&d ~=~  \half   ~=~d_{23}  \ , \quad a_1 = a_{2} = 0 \ , \nn\\
&&d ~=~  0   ~=~ d_{3}  \ , \quad a_1 = a_{2} = 0 \   .\nn \eea
For ~$d>d_1$~ there are no singular vectors and we have unitarity.
At ~$d=d_1$~ there is a singular vector:
\bb  v^1_{\d_1} ~&=&~ \Big( a_1  (a_1+a_2+1) X_{\d_1} ~-~ a_1  X_{\d_3}\ X_{13} ~-~   (a_1+a_2+1) X_{\d_2}\ X^+_{1} ~+\nn\\ &&+~  X_{\d_3}\ X^+_{2}\ X^+_{1} \Big)\,v_0  \label{svn3d1}
\eea
which is given in PBW basis, where ~$X_{\d_j} \in \cg^+$~ are the vectors corresponding to the weight vectors
~$\d_j\,$, ~ $X_{13}$~ is the compact root vector for ~$\a_{13} = \a_1+\a_2 = \d_1-\d_3\,$.
This singular vector is non-trivial for ~$a_1\neq 0$~ and must be eliminated to obtain an UIR.
Below ~$d<d_1$~ there is no unitarity for ~$a_1\neq 0$. On the other hand for $a_1= 0$ the singular vector $v^1_{\d_1}$ is
descendant of the compact root singular vector  ~$X_1^+\,v_0$~ which is already factored out for ~$a_1 = 0$.
Thus, below we discuss only the cases with ~$a_1  = 0$.\nl
The singular vector at ~$d ~=~ d_{12}$~ corresponding to the root ~$\d_1 + \d_2 = \a_1 + 2\a_2 + 2\a_3$~ is:
\eqnn{vs12} v^{1}_{\d_1 + \d_2} &=& \bigg( X_{\delta _3} X_{\delta _2} X^+_2X^+_1+
\frac{1}{2}\, (X_{\delta _3})^2  (X^+_2)^2 X^+_1 -a_1\, (X_{\delta _3})^2  X^+_2 X_{13}
   + \nn\\ &&+\ 2\left(a_2+1\right)\, X_{\delta _3} X_{\delta _2} X_{13} -2\left(a_1+a_2+1\right)\, X_{\delta _3} X_{\delta _1} X^+_2  + \nn\\ &&+\ \left(a_1+1\right) \left(a_1+a_2+1\right)\, X_{\delta _1+\delta _3} X^+_2 +4a_2 \left(a_1+a_2+1\right) \, X_{\delta _2} X_{\delta _1} + \nn\\ &&\ + 2a_2
   \left(a_1+a_2+1\right)\, ( X_{\delta _2})^2 X^+_1  -\frac{1}{2} \left(a_1+2
   a_2+1\right)\, X_{\delta _2+\delta _3} X^+_2X^+_1 -\nn\\ && -\  \left(-a_2 a_1+a_1+a_2+1\right)\, X_{\delta _2+\delta _3} X_{13}- \nn \\ && -\ 2 \left(a_1+1\right) a_2
   \left(a_1+a_2+1\right) X_{\delta _1+\delta _2} \bigg)v_0  \eea
with norm:
$$16 \left(2 d-a_2-3\right) \left(a_1^2+2 a_1+2
d-a_2-2\right) a_2 \left(a_2+1\right) \left(a_1+a_2+1\right)
\left(a_1+a_2+2\right).$$
For ~$a_1 = 0,\  a_2\neq 0$~ it is non-trivial and gives rise to a invariant subspace which must be factored out for unitarity.
For ~$d<d_{12}$~ the vector \eqref{vs12}  has negative norm and there is
no unitarity for ~$a_2\neq 0$, except at the isolated unitary point ~$d ~=~ 1 + \half a_2 = d_2 > d_{13}$.
At this point there is a singular vector ~$v_{s2}$, while the vector \eqref{vs12}
 is descendant of compact root singular vector ~$X_1\,v_0\,$~ and ~$v_{s2}\,$.\nl
Further, we consider ~$a_1 = a_2= 0$. Then   the vector
~$v^1_{\d_1+\d_2}$~ is descendant of compact root singular vectors
~$X^+_1\,v_0\,$~ and ~$X^+_2\,v_0\,$, thus, there is no obstacle for
unitarity for ~$1 <d$. The next reducibility points (coinciding here) are  ~$d=d_{13}=d_2=1$.
The singular vector for ~$d=d_2$~ and ~$m=1$~
has weight ~$\d_2=\a_2+\a_3$~ and is given by:
\eqn{svd22} v^1_{\d_2} ~=~ \Big( a_2 X^+_2 X^+_3 - (a_2+1) X^+_3 X^+_2\Big) v_0 \ee
For ~$a_2=0$~ it is a descendant of the compact root singular vector ~$X^+_2\, v_0\,$.
The singular vector for ~$d=d_{13}=1 $~ and ~$m=1$~ has weight ~$\d_1+\d_3=a_1+\a_2+2\a_3$~ \cite{Dos}~:
\eqnn{svd133} v^1_{\d_1+\d_3} ~&=&~ \Big(ha_1\, X^+_1 (X^+_3)^2 X^+_2 +a_1\,X^+_1 X^+_3 X^+_2 X^+_3 -\nn\\
&&-\ h(a_1+1) (X^+_3)^2 X^+_2 X^+_1\ - \nn\\
&&-\ ha_1\, X^+_1  X^+_2 (X^+_3)^2 -  (a_1+1) X^+_3 X^+_2 X^+_3 X^+_1\ +\nn\\
&&+\ h(a_1+1) X^+_2 (X^+_3)^2  X^+_1 \Big) ,   \\
 && \qquad \qquad h= 1+\half ( a_1+a_2) \nn \eea
 The above vector is given in the simple root basis most appropriate for the case.
 For ~$a_1=0$~ it is a descendant of the compact root singular vector ~$X^+_1\, v_0\,$. However, there is a subsingular vector:
\eqn{vssd1} v_{ss} ~=~ (X_{\d_1}X_{\d_2}X_{\d_3}-X_{\d_3}X_{\d_2}X_{\d_1})v_0 \ee
with norm: ~$16d(d-1)(2d-1)$. This must be factorized in order to obtain UR.  Then for ~$\half < d <1$~ there will be no unitarity due to the last norm.\\
\indent Finally, at the next reducibility point: ~$d=d_{23}=\half$~ there is a singular vector of weight ~$\d_2+\d_3=\a_2 + 2\a_3$~:
\eqn{svd23} v^1_{\d_2+\d_3} ~=~ (2X_{\d_2+\d_3} - 4 X_{\d_2}X_{\d_3} + X_{2\d_3}X^+_2 )v_0 \ee
It should be factored out to get unitarity. The subsingular vector \eqref{vssd1} has zero norm for ~$d=\half$~ and furthermore
  it is descendant of ~$v^1_{\d_2+\d_3}$~ and the compact root singular vector ~$X^+_2\, v_0\,$. Finally, for ~$d<\half$~ there is no unitarity since then the norm of \eqref{svd23} is negative, except at the trivial isolated unitary point
~$d=0=a_1 = a_2$~ of one-dimensional irrep.\bsq

\section{Character formulae}

\subsection{Character formulae: generalities}

\nt
In the beginning of this subsection we follow \cite{Dix}.
Let ~$\hcg$~ be a simple Lie algebra of rank ~$\ell$~
with Cartan subalgebra
~$\hch$, root system ~$\hD$, simple root system ~$\hpi$.
Let ~$\G$, (resp. $\G_+$), be the set of all integral, (resp.
integral dominant), elements of $\hch^*$, i.e., $\l
\in \hch^*$ such that $(\l , \a_i^\vee) \in \bbz$, (resp. $\bbz_+$),
for all simple roots $\a_i\,$, ($\a_i^\vee \equiv 2\a_i/(\a_i,\a_i)$).
Let $V$ be a lowest weight module with lowest weight $\L$ and
lowest weight vector $v_0\,$. It has the following decomposition:
\eqn{wei} V ~=~ \mathop{\oplus}\limits_{\mu\in\G_+} V_\mu ~~, \ \
~~~V_\mu ~=~ \{ u \in V ~\vert ~Hu = (\l + \mu)(H)u, \
\forall ~H\in\ch \} \ee
(Note that $V_0 = \bbc v_0\,$.) Let $E(\ch^*)$ be
the associative abelian algebra consisting of the series
$\sum_{\mu \in \ch^*} c_{\mu} e(\mu)$ , where $c_{\mu} \in \bbc ,
 ~c_{\mu} = 0$ for $\mu$ outside the union of a finite number of
sets of the form $D(\l) = \{ \mu \in \ch^* \vert \mu \geq
\l \}$~, ~using some ordering of $\ch^*$,
e.g., the lexicographic one; the formal
exponents $e(\mu)$ have the properties:~
$e(0) = 1, \ e(\mu) e(\nu) = e(\mu + \nu)$.

Then the (formal) character of $V$ is defined by:
\eqn{cha}ch_0~V ~~=~~ \sum_{\mu \in\G_+} (\dim \ V_\mu) ~e(\L+\mu)
 ~~=~~ e(\L) \sum_{\mu\in\G_+} (\dim \ V_\mu) ~e(\mu) \ee
(We shall use subscript '0' for the even case.)

For a Verma module, i.e.,
~$V = V^\L$~ one has $\dim \ V_\mu = P(\mu)$,
where ~$P(\mu)$ is a generalized
partition function, $P(\mu) = \#$ of ways $\mu$ can be presented
as a sum of positive roots $\beta$, each root taken with its
multiplicity $\dim \cg_{\beta}$ ($=1$ here),
$P(0)\equiv 1$. Thus, the character formula for Verma modules is:
\eqn{chv}ch_0~V^\L ~~=~~ e(\L)\sum_{\mu\in\G_+} P(\mu) e(\mu) ~~=~~ e(\L)
\prod_{\a \in\D^+}(1 - e(\a))^{-1}. \ee

Further we recall the
standard reflections in $\hch^*$~:
\eqn{rfl} s_\a(\l) ~=~ \l - (\l, \a^\vee)\a \,,
\quad \l \in \hch^* \,, \quad \a\in\hD. \ee
The Weyl group ~$W$~ is generated by the simple reflections
$s_i \equiv s_{\a_i}$, $\a_i\in\hpi\,$.
Thus every element $w\in W$ can
be written as the product of simple reflections. It is said that
$w$ is written in a reduced form if it is written with the
minimal possible number of simple reflections; the number of
reflections of a reduced form of $w$ is called the length of $w$,
denoted by $\ell (w)$.

The Weyl character
formula for the finite-dimensional irreducible LWM $L_\L$
over ~$\hcg$, i.e., when ~$\L\in -\G_+\,$,
has the form:
\eqn{chm}
ch_0~L_\L ~~=~~ \sum_{w\in W}(-1)^{\ell(w)} ~ch_0~V^{w\cdot\L} \,, \quad
\L\in -\G_+ \ee
where the dot ~$\cdot$~ action is defined by $w\cdot \l = w(\l -
\r) +\r$.
For future reference we note:
\eqn{rft} s_\a ~\cdot ~\L ~~=~~ \L ~+~ n_\a \a  \ee
where
\eqn{nchs}  n_\a ~~=~~ n_\a(\L) ~~\doteq~~ (\r-\L,\a^\vee)
~~=~~ (\rho-\L)(H_\a)
\,, \quad \a\in\D^+ .  \ee

In the case of basic classical Lie superalgebras the first
character formulae were given by Kac \cite{Kc,Kacc}.\footnote{Kac considers
highest weight modules but his results are immediately
transferable to lowest weight modules.} For all such
superalgebras  -- except $osp(1/2n)$ --
the character formula for Verma modules is \cite{Kc,Kacc}:
\eqn{chvv} ch~V^\L
~=~ e(\L)
\( \prod_{\a \in\D^+_\o}(1 - e(\a))^{-1} \)
\( \prod_{\a \in\D^+_\I}(1 + e(\a)) \).
\ee
We are however interested exactly in the ~$osp(1/2n)$~ when the Verma module character
formula is:
\eqn{chvvo} ch~V^\L
~=~ e(\L)
\( \prod_{\a \in\bd^+}(1 - e(\a))^{-1} \)
\ee
Naturally, the character
formula for the finite-dimensional irreducible LWM $L_\L$ is again
\eqref{chm} using the Weyl group ~$W_n$~ of ~$B_n\,$.

\subsection{Multiplets}

A Verma module $V^\L$ may be reducible w.r.t. to many positive roots, and thus there maybe many
Verma modules isomorphic to its submodules. They themselves may be reducible, and so on.

One main ingredient of the  approach of \cite{Dob} is as follows. We group the
(reducible) Verma modules with the same Casimirs in sets called ~{\it
multiplets} \cite{Dobmul}. The multiplet corresponding to fixed
values of the Casimirs may be depicted as a connected graph, the
vertices of which correspond to the reducible Verma modules and the lines
between the vertices correspond to embeddings between them. The
explicit parametrization of the multiplets and of their Verma modules is
important for understanding of the situation.

If a Verma module $V^\L$ is reducible w.r.t. to all simple roots (and thus w.r.t. all positive roots),
i.e., ~$m_k\in\bbn$~ for all $k$, then the irreducible submodules are isomorphic to the
  finite-dimensional irreps of $\cg^\bac$ \cite{Dob}.
(Actually, this is a condition only for $m_n$ since ~$m_k\in\bbn$~ for   $k=1,\ldots,n-1$.)
  In these cases we  have the ~{\it main
  multiplets}~ which are isomorphic to the Weyl group of $\cg^\bac$ \cite{Dob}.

In the cases of non-dominant weight ~$\L$~ the character formula for
the irreducible LWM is \cite{KL}~:
\eqn{chg}
ch~L_\L ~~=~~ \sum_{w\in W \atop w\leq w_\L}
~(-1)^{\ell(w_\L w)} ~P_{w,w_\L}(1) ~ch~V^{w\cdot (w_\L^{-1}\cdot\L)}
\,, \quad \L\in \G \ee
where ~$P_{y,w}(u)$~ are the Kazhdan--Lusztig polynomials
$y,w\in W$ \cite{KL} (for an easier exposition see \cite{Do-KL}),
~$w_\L$~ is a unique element of $W$ with minimal length
such that the signature of ~$\L_0 ~=~ w^{-1}_\L \cdot \L$~
is anti-dominant or semi-anti-dominant:
\eqn{sngg} \chi_0 ~~=~~ (m'_1,\ldots,m'_n),
\qquad m'_k ~=~ 1 - \L_0(H_k) ~\in \bbz_- \ . \ee
Note that ~$P_{y,w}(1) \in\bbn\,$ for $y\leq w$.

When ~$\L_0$~ is semi-anti-dominant, i.e., at least one ~$m'_k=0$,
then in fact ~$W$~ is replaced by a reduced Weyl group ~$W_R\,$.

Most often the value of ~$P_{y,w}(1)$~ is equal to 1 (as in the character
formula for the finite-dimensional irreps), while the cases ~$P_{y,w}(1)>1$~
are related to the appearance of subsingular vectors, though the situation
is more subtle, see \cite{Do-KL}.

It is interesting to see how the reducible points relevant for unitarity fit in the multiplets.
In the case of ~$d_{ij}$~ \eqref{boun} and using  \eqref{cftw} we have:
\eqn{bounm}
m_n (d_{ij}) ~=~ 1 - 2m_j - \cdots -2m_{n-1} - m_i - \cdots - m_{j-1} \ .\ee

In the case of ~$d_{i}$~ \eqref{bounz}   we have:
\eqn{bounzm}
m_n (d_i) ~=~ 1 - 2m_i - \cdots - 2m_{n-1} \ . \ee

As expected the weights related to positive energy $d$ are not dominant
($m_n (d_{ij})\in\bbz_-$, ~$m_n (d_{i})\in -\bbn$, ($i<n$)), since the positive energy UIRs are
infinite-dimensional. (Naturally, $m_n (d_{n})=1$ falls out of the picture since $d_n<0$.)

Thus, the Verma modules with weights related to positive energy would be somewhere in the
main multiplet (or in a reduction of the main multiplet), and the first task for calculating
the character is to find the ~$w_\L$~ in the character formula \eqref{chg}. This we do in the next
subsection in the case ~$n=3$.


\subsection{The case n=3}

In order to illustrate what the main ideas we consider the first non-trivial example ~$n=3$, i.e.,
$osp(1/6)$ actually using ~$\boldmath{B_3}$. The Weyl group $W_n$ of $B_n$ has ~$2^n n!$~ elements,
i.e., 48 for $B_3$.  Let ~$S ~=~ (s_1,s_2,s_3)$, ~$s_i \equiv s_{\alpha_i}\,$, be the simple
reflections.  They fulfill the
following relations:
\eqn{defrel} s_1^2 = s_2^2 =s_3^2 = e, ~(s_1s_2)^3 = e, ~(s_2s_3)^4 = e,
~s_1s_3 ~=~ s_3s_1\,,\ee $e$ ~being the identity of ~$W_3\,$.
The 48 elements may be listed as:
\eqnn{weylel} && e\,, ~s_1\,, ~s_2\,, ~s_3\, \\ &&
s_1s_2\,, ~ s_1s_3\,, ~s_2s_1\,, ~s_2s_3\,, ~s_3s_2\,, ~ \nn\\ &&
 s_1s_2s_1\,, ~s_1s_2s_3\,, ~s_1s_3s_2\,, ~s_2s_1s_3\,, ~ s_2s_3s_2\,, ~s_3s_2s_1\,,   ~s_3s_2s_3\,,
 \nn\\ &&
s_1s_2s_1s_3\,, ~s_1s_2s_3s_2\,, ~s_1s_3s_2s_1\,, ~ s_1s_3s_2s_3\,, ~\nn\\ &&
s_2s_3s_2s_1\,, ~s_2s_1s_3s_2\,, ~s_3s_2s_3s_1\,,~ s_3s_2s_3s_2\,,~\nn\\ &&
s_1s_2s_3s_2s_1\,, ~s_1s_3s_2s_1s_3\,, ~s_1s_2s_1s_3s_2\,, ~ s_1s_3s_2s_3s_2\,, ~\nn\\ &&
s_2s_1s_3s_2s_1\,, ~s_2s_1s_3s_2s_3\,, ~ s_3s_2s_3s_1s_2\,,~ s_3s_2s_3s_2s_1\,,~\nn\\ &&
s_1s_3s_2s_3s_2s_1\,, ~s_1s_3s_2s_1s_3s_2\,, ~ s_1s_2s_1s_3s_2s_1\,, ~s_2s_1s_3s_2s_1s_3\,, ~\nn\\ &&
s_2s_1s_3s_2s_3s_2\,, ~ s_3s_2s_3s_1s_2s_1\,,~ s_3s_2s_3s_1s_2s_3\,,~\nn\\ &&
s_2s_1s_3s_2s_3s_2s_1\,, ~ s_2s_1s_3s_2s_3s_1s_2\,, ~ s_3s_2s_1s_2s_3s_2s_1\,,~\nn\\ &&
s_3s_2s_3s_1s_2s_1s_3\,,~ s_3s_2s_3s_1s_2s_3s_2\,,~~\nn\\ &&
s_2s_3s_2s_1s_2s_3s_2s_1\,,~ s_3s_2s_1s_3s_2s_3s_2s_1\,, ~s_3s_2s_1s_3s_2s_3s_1s_2\,, ~\nn\\ &&
s_2s_3s_2s_1s_3s_2s_3s_2s_1\ .
\nn\eea
\fig{}{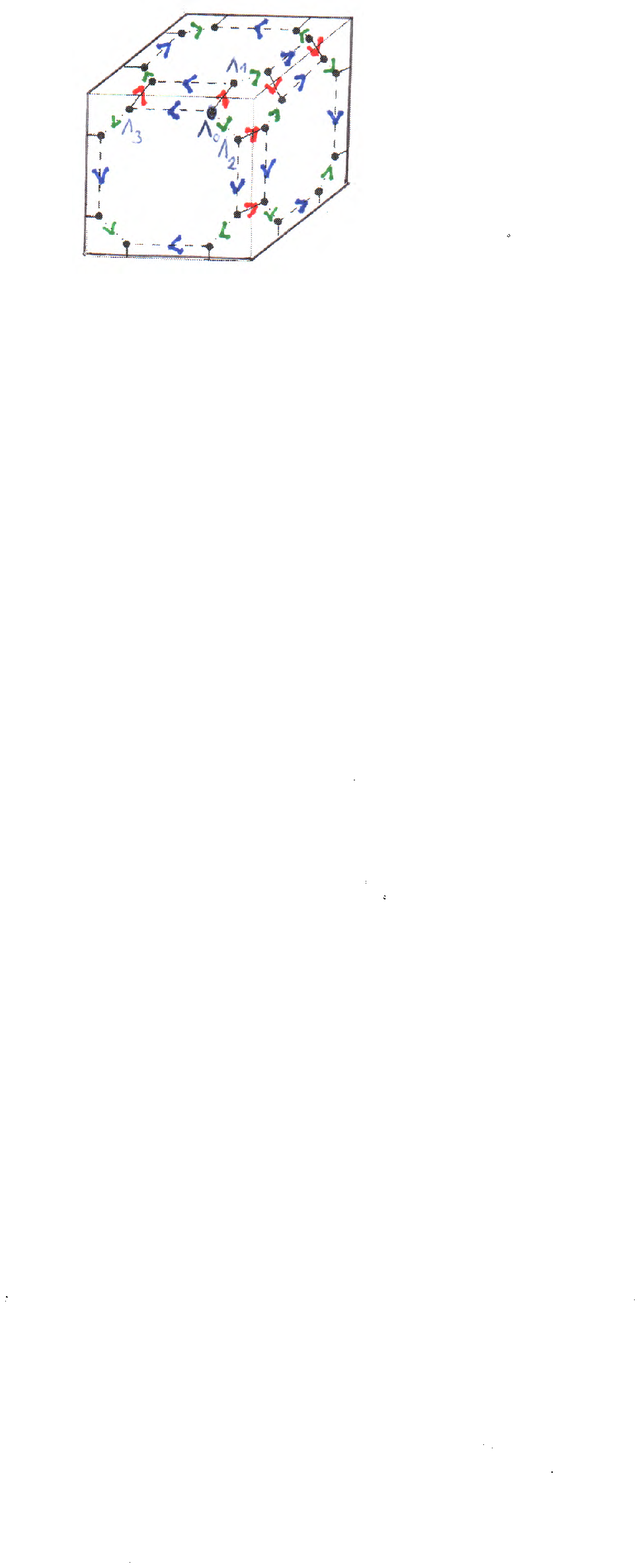}{13cm}
This Weyl group may be pictorially represented on a cube as in the figure,
where we have given only the simple root reflections, namely,
continuous (red) arrows represent action of reflection $s_1$, dashed (blue) arrows represent action of reflection $s_2$, dotted (green) arrows represent action of reflection $s_3$.
Each face of the cube contains eight elements related by blue and green arrows representing the Weyl group of ~$B_2$~ generated by $s_2$ and $s_3\,$. The figure contains also eight sextets (around the eight corners of the cube). Each sextet is related by red and green arrows representing the Weyl group of ~$A_2$~ generated by $s_1$ and $s_2\,$. Finally there are 12 quartets (straddling the edges of the cube).  Each quartet is formed by red and blue arrows representing the Weyl group of ~$A_1\times A_1$~ generated by the commuting reflections $s_1$ and $s_3\,$.

We use the same diagram to depict the main multiplets containing the Verma modules $V^{\L_0}$ which contain
(as factor module) the finite-dimensional irreps of $B_3$, i.e.,
with dominant weights $\L_0$, i.e., with Dynkin labels\\ $(m_1,m_2,m_3)$, ~$m_k\in\bbn$.
We may do this since these multiplets are isomorphic to the Weyl group, $W_3$ in our case.
On the picture we have indicated the modules, $\L_0$ and $\L_k ~=~ s_k \cdot \L_0\,\,$, $k=1,2,3$. The mentioned isomorphism is fixed by assigning to $\L_0$ the identity element ~$e$~ of $W_3$, and to ~$\L_k$~ the reflections ~$s_k\,$.

The character formula for the Verma modules in our case is given explicitly by:
\eqnn{chvvoz} ch~V^\L
&=& \frac{e(\L)}{ (1-t_1)(1-t_2) (1-t_1t_2)}\ \times \\
&\times& \frac{1}{
 (1-t_3) (1-t_2t_3)(1-t_1t_2t_3) (1-t_2t^2_3) (1-t_1t_2t^2_3) (1-t_1t^2_2t^2_3)}
\nn\eea where ~$t_j ~\equiv~ e(\a_j)$.

Now we give the character formulae  of the
five boundary or isolated unitarity cases.
Below we shall denote the signature of the dominant weight ~$\L_0$~ which determines
the main multiplet by ~$(m'_1,m'_2,m'_3)$, ~$m'_k\in\bbn$,
 using primes to distinguish from the signatures
of the weights we are interested. We shall use also reductions of the main multiplet
when the weights are semi-dominant, i.e., when some ~$m'_k=0$.

\vsh


\nt\bu In the case of ~$d=d_1=2+\ha(a_1+a_2)$~ there are twelve members of the multiplet which
is a submiltiplet of a main multiplet.
(Remember that that ~$m_1>1$~ since $a_1\neq 0$.)
They are grouped into two standard ~$sl(3)$~ submultiplets
of six members. The first submultiplet starts from ~$V^{\L_0^{d_1}}$,
where ~$\L_0^{d_1} ~=~ w\cdot \L_0$, ~$w ~=~ w_{\L_0^{d_1}} ~=~
s_{2}s_{1}s_{3}s_{2}s_{3}\,$,
with signature:
\eqn{weightsone}
  \L_0^{d_1} ~:~ (m_1,m_2,m'_3=1-2m_{12})\ , \quad m_1,m_2\in\bbn \ , \quad m_{12} \equiv m_1+m_2\
\ .   \ee
The other submultiplet starts from ~$V^{\L'_0}$~ with   ~$\L'_0 ~=~ \L_0^{d_1} + \d_1 ~=~
\L_0^{d_1} + \a_1+ \a_2+ \a_3$,
with signature: ~$\L'_0 ~:~ (m_1-1,m_2,m'_3=1-2m_{12})$, $m_1>1$.
The character formula is  ~\eqref{chg}~ with ~$w_\L ~=~ w_{\L_0^{d_1}}$~:
\eqnn{chardone} && {\rm ch}\, \L_0^{d_1} =
\frac{e(\L_0^{d_1})}{ (1-t_3) (1-t_2t_3)(1-t_1t_2t_3) (1-t_2t^2_3) (1-t_1t_2t^2_3) (1-t_1t^2_2t^2_3)}  \nn\\ &&\nn\\
&&\times\ \{ \ {\rm ch}\, \L_{m_1,m_2}(t_1,t_2) ~-~ t_1t_2t_3\, {\rm ch}\, \L_{m_1-1,m_2}(t_1,t_2) \ \}\ , ~~m_1>1 \eea
where ~${\rm ch}\, \L_{m_1,m_2}(t_1,t_2)$~ is the normalized  character of the finite-dimensional ~$sl(3)$~ irrep with Dynkin
labels ~$(m_1,m_2)$~ (and dimension ~$m_1m_2(m_1+m_2)/2$):
\eqn{charfd} {\rm ch}\ \L_{m_1,m_2}(t_1,t_2) ~=~
\frac{1 - t_1^{m_1} - t_2^{m_2} + t_1^{m_1} t_2^{m_{12}} + t_1^{m_{12}} t_2^{m_2}  - t_1^{m_{12}} t_2^{m_{12}}}{(1-t_1) (1-t_2)(1-t_1t_2)} \ee
Naturally, the latter formula is a polynomial in ~$t_1,t_2\,$, e.g., ~$ {\rm ch}\ \L_{1,1}(t_1,t_2) ~=~ 1$. Note that \eqref{chardone}
trivializes for ~$m_1=1$~ since the second term disappears by the formal substitution: ~${\rm ch}\ \L_{0,m_2}(t_1,t_2) ~=~ 0$.

 \vsh

\nt\bu In the case of ~$d=d_{12}=\ha(3+a_2)$~ which is relevant for unitarity, i.e., $m_1 =  1$,
  there are again twelve  members of the multiplet. The corresponding  signature is:
\eqn{weightsonetwo}
 \L_0^{d_{12}} ~:~ (1,m_2,m'_3=-2m_{2})\ , \quad m_2\in\bbn
\ .    \ee
  The multiplet is submultiplet
  of a reduced multiplet with 24 members obtained from a main multiplet for ~$m'_3=0$.
As above our multiplet consists of  two standard ~$sl(3)$~ submultiplets
of six members. The first submultiplet starts from ~$V^{\L_0^{d_{12}}}$,
where ~$\L_0^{d_{12}} ~=~ w\cdot \L_0$, ~$w ~=~ w_{\L_0^{d_{12}}} ~=~
s_{3}s_{2}s_{1}\,$.
The other submultiplet starts from ~$V^{\L'_0}$~ with
~$\L'_0 ~=~ \L_0^{d_{12}} + m_2(\a_1 + 2\a_2 + 2\a_3) ~=~ \L_0^{d_{12}} + m_2(\d_1 + \d_2)$~
with signature: ~$\L'_0 ~:~ (1,m_2-1,-2m_{2})$.
The character formula is ~\eqref{chg},   with ~$W ~\mt W_R\,$, (where ~$W_R$~ is  a
reduced 24-member Weyl group) and with ~$w_\L ~=~w_{\L_0^{d_{12}}}$~:
\eqnn{chardonetwo} && {\rm ch}\, \L_0^{d_{12}} =
\frac{e(\L_0^{d_{12}})}{ (1-t_3) (1-t_2t_3)(1-t_1t_2t_3) (1-t_2t^2_3) (1-t_1t_2t^2_3)
(1-t_1t^2_2t^2_3)} \nn\\ &&\nn\\
&&\times \{ \ {\rm ch}\, \L_{1,m_2}(t_1,t_2) ~-~ (t_1 t^2_2 t_3^2)^{m_2}\,
 {\rm ch}\, \L_{1,m_2-1}(t_1,t_2) \ \}\ , ~~m_2>1 \eea
where ~${\rm ch}\, \L_{m_1,m_2}$~ are the $sl(3)$ characters defined in \eqref{charfd}.

\vsh

\nt\bu In the case of ~$d=d_2=1+\ha a_2\geq d_{13}\,$, i.e., $m'_3=1-2m_2$,
the corresponding  signature is:
\eqn{weightstwoo}
 \L_0^{d_2}  ~:~ (m_1,m_2,m'_3=1-2m_2) \ . \ee
We should  consider two subcases
$$1+m_{1}-m_2 > 0  \qquad  {\rm or}  \qquad  1+m_{1}-m_2 \leq 0$$

We start with the ~{\it first subcase}~ which is relevant when ~$d=d_2=d_{13}=1$~ and ~$a_1=a_2=0$,
then   ~$m_1=m_2=1$, and the signature is:
\eqn{weightstwo}
 \L_0^{d_2=d_{13}}  ~:~   (1,1,-1) \  . \ee

Our multiplet is
a submultiplet of a 12-member reduced multiplet obtained when the signature of ~$\L_0$~ is
~$(m'_1,m'_2,m'_3)  ~=~
(1,0,1)$, and then $\L_0^{d_2=d_{13}}$ is a submodule with signature \eqref{weightstwo}.
Thus, we have ~$\L_0^{d_{2}=d_{13}}  ~=~  s_3    \cdot \L_0\,$, i.e.,
~$w_{\L_0^{d_{2}=d_{13}}} ~=~ s_3\,$.

Explicitly, our 12-member multiplet has two $sl(3)$ submultiplets. First we take into account a  $sl(3)$
sextet starting from ~$\L_0^{d_2=d_{13}}$~ with parameters ~$(1,1)$. Then
there is   a  $sl(3)$ sextet starting from ~$\L_0^{d_2=d_{13}}  + \a_1 + 2\a_2 + 3\a_3$~
with parameters ~$(1,1)$.  Note that that ~$\a_1 + 2\a_2 + 3\a_3 = \d_{1} + \d_{2} + \d_{3}$~
is the weight of the subsingular vector \eqref{vssd1}.


The character formula is ~\eqref{chg},   with ~$W ~\mt W_R\,$, (where ~$W_R$~ is  a
reduced 12-member Weyl group) and ~$w_\L = s_3$~:
\eqnn{chardonezz} &&  {\rm ch}\, \L_0^{d_2=d_{13}} ~=~ \nn\\
&&=~ \frac{e(\L_0^{d_2=d_{13}})} { (1-t_3) (1-t_2t_3)(1-t_1t_2t_3) (1-t_2t^2_3) (1-t_1t_2t^2_3)
(1-t_1t^2_2t^2_3)}\times    \nn\\ &&\nn\\
&&\times\ \{ \  1  
~-~ t_1t^2_2t_3^3 \ \}  
\eea

\vsh

\nt\bu In the case of ~$d=d_{2}=1+\half a_2 > d_{13}=1$, i.e.,
~$m_1=1$, ~$m_2=1+  a_2>1$, thus, this is the subcase
~$1+m_{1}-m_2 = m_{13} \leq 0$. The  multiplet has 24 members for ~$m_2 >2$~
($m_{13}<0$)
and starts with ~$\L^{'d_2}_0 ~=~ ~=~ s_3 \, s_2 \, s_1\cdot \L_0\,$,
 with signatures:
\eqnn{weightstwoz}
&& \L_0 ~:~ (m_2-2,1,1) \ , \nn\\
&& \L^{'d_2}_0  ~:~ (1,m_2,m'_3=1-2m_2)\ , \quad m_2\in 1+\bbn \   . \eea
 It has four $sl(3)$ submultiplets. First we take into account a  $sl(3)$
sextet starting from ~$\L^{'d_2}_0$~ with parameters ~$(1,m_2)$. Then
there is   a  $sl(3)$ sextet starting from ~$\L^{'d_2}_0 + \a_{23}$~ 
with parameters ~$(2,m_2-1)$.
Then there  is a $sl(3)$ sextet with
parameters ~$(2,m_2-2)$~ starting from a Verma module ~$V^{\L''}$, 
~$\L'' = \L^{'d_2}_0 +  \a_1 + 3  \a_{23} $.
Finally,   there is a    $sl(3)$ sextet with parameters ~$(1,m_2-2)$,
starting from a Verma module ~$V^{\L'''}$, 
  ~$\L''' = \L^{'d_2}_0 +  2(\a_1+2\a_2+2\a_3)$.\nl
  We have the  ~{\it Conjecture}~ that  the
  character formula is ~\eqref{chg}~ and ~$w_\L = s_3 \, s_2 \, s_1$~:
\eqnn{chh12}  &&  {\rm ch}\, \L^{'d_2}_0 = \frac{e(\L^{'d_2}_0)}{
(1-t_3) (1-t_2t_3)(1-t_1t_2t_3) (1-t_2t^2_3) (1-t_1t_2t^2_3)
(1-t_1t^2_2t^2_3)}    \nn\\ &&\nn\\
&&\qquad \times\ \{ \ {\rm ch}\, \L_{1,m_2}(t_1,t_2) ~-~
  t_2 t_3\, {\rm ch}\, \L_{2,m_2-1}(t_1,t_2) ~+ \nn\\ && \nn\\ &&
\qquad +~ t_1 t_2^3 t_3^3 \, {\rm ch}\, \L_{2,m_2-2}(t_1,t_2)
~-~ t_1^2 t_2^4 t_3^4 \, {\rm ch}\, \L_{1,m_2-2}(t_1,t_2)\ \} \eea
When ~$m_2=2$~ ($a_2=1$,~ $m_{13}=0$) the weight ~$\L_0$~ is semi dominant, the main multiplet reduces
to 24 members, our multiplet reduces to only 12 members,
consisting of the first two $sl(3)$ submultiplets mentioned above.
The   character formula takes this into account by construction
since for $m_2=2$ the terms in the 2nd row are automatically zero
(due to the fact that the $sl(3)$ character formula gives zero:
~${\rm ch}\, \L_{1,0}(t_1,t_2) ~=~0$).

\vsh

\nt\bu
In the case of ~$d=d_{23}=\half$, ~$a_1=a_2=0$, i.e.,  ~$m_1=m_2=1$,
and the signature is:
\eqn{twothree} \L_0^{d_{23}} ~:~ (1,1,0) \ .\ee  
This is in fact a multiplet with 24 members which is reduction of the main multiplet
starting with the semi dominant weight \eqref{twothree}.

The multiplet consists of  four $sl(3)$ submultiplets. First  there is a  $sl(3)$
sextet starting from ~$\L^{d_{23}}_0$~ with parameters ~$(1,1)$. Then
a  $sl(3)$ sextet starting from ~$\L^{d_{23}}_0~+~ \a_2 + 2\a_3$~ with parameters ~$(2,1)$.
 Then a  $sl(3)$  
sextet starting from ~$\L^{d_{23}}_0~+~ \a_1 + 2\a_2 + 4\a_3$~ with parameters ~$(1,2)$.
 Then a  $sl(3)$ 
sextet starting from ~$\L^{d_{23}}_0~+~ 2\a_1 + 4\a_2 + 6\a_3$~ with parameters ~$(1,1)$.

The character formula is ~\eqref{chg}, however, with ~$W ~\mt W_R\,$,~ where ~$W_R$~ is the
reduced 24-member Weyl group, (generated by ~$s_1,s_2,s_3s_2s_3\,$) and ~$w_\L  =1$~:
\eqnn{chh12z}  &&  {\rm ch}\, \L^{d_{23}}_0 ~=~ \frac{e(\L^{d_{23}}_0)}{
(1-t_3) (1-t_2t_3)(1-t_1t_2t_3) (1-t_2t^2_3) (1-t_1t_2t^2_3)
(1-t_1t^2_2t^2_3)}    \nn\\ &&\nn\\
&&\qquad\qquad \times\ \{ \  1 ~-~t_2 t_3^2 \, {\rm ch}\,
\L_{2,1}(t_1,t_2)
   ~+ \nn\\ && \nn\\ &&
\qquad\qquad ~+~ t_1 t^2_2 t_3^4 \, {\rm ch}\, \L_{1,2}(t_1,t_2) ~-~
t_1^2 t_2^4 t_3^6  \ \} ~= \nn\\ && \nn\\ && ~=~
\frac{e(\L^{d_{23}}_0)}{ (1-t_3) (1-t_2t_3)(1-t_1t_2t_3)
(1-t_2t^2_3) (1-t_1t_2t^2_3)
(1-t_1t^2_2t^2_3)}    \nn\\ &&\nn\\
&&\qquad\qquad \times\ \{ \  1 ~-~ t_2 t_3^2 \, (1+t_1+t_1t_2)
     ~+ \nn\\ && \nn\\ &&
\qquad\qquad ~+~  t_1 t^2_2 t_3^4 \,(1+t_2+t_1t_2) ~-~ t_1^2 t_2^4
t_3^6  \ \} ~=\nn\\ && \nn\\ && =~
\frac{e(\L^{d_{23}}_0)}{
(1-t_3) (1-t_2t_3)(1-t_1t_2t_3) }
\eea

\vspace{5mm}

\section*{Acknowledgements}

V.K. Dobrev is supported in part by Bulgarian NSF Grant DFNI T02/6. I.~Salom is supported in part by the Serbian Ministry of Science and Technological Development under grant number OI 171031.

\vspace{5mm} 

\end{document}